\documentstyle[12pt, amsfonts]{amsart}
\begin{document}
\def\R{{\mathbb R}}
\def\Z{{\mathbb Z}}
\def\C{{\mathbb C}}
\newcommand{\trace}{\rm trace}
\newcommand{\Ex}{{\mathbb{E}}}
\newcommand{\Prob}{{\mathbb{P}}}
\newcommand{\E}{{\cal E}}
\newcommand{\F}{{\cal F}}
\newtheorem{df}{Definition}
\newtheorem{theorem}{Theorem}
\newtheorem{lemma}{Lemma}
\newtheorem{pr}{Proposition} 
\newtheorem{co}{Corollary}
\newtheorem{problem}{Problem}
\def\n{\nu}
\def\sign{\mbox{ sign }}
\def\a{\alpha}
\def\N{{\mathbb N}}
\def\A{{\cal A}}
\def\L{{\cal L}}
\def\X{{\cal X}}
\def\F{{\cal F}}
\def\c{\bar{c}}
\def\v{\nu}
\def\d{\delta}
\def\diam{\mbox{\rm dim}}
\def\vol{\mbox{\rm Vol}}
\def\b{\beta}
\def\t{\theta}
\def\l{\lambda}
\def\e{\varepsilon}
\def\colon{{:}\;}
\def\pf{\noindent {\bf Proof :  \  }}
\def\endpf{ \begin{flushright}
$ \Box $ \\
\end{flushright}}

\title[slicing inequality]
{A remark on measures of sections of $L_p$-balls}

\author{Alexander Koldobsky}

\address{Department of Mathematics\\ 
University of Missouri\\
Columbia, MO 65211}

\email{koldobskiya@@missouri.edu}

\thanks{The first named author was partially supported by the US National Science Foundation, 
grant DMS-1265155.}

\author{Alain Pajor}

\address{Universit\'e Paris-Est,
Laboratoire d'Analyse et Math\'ematiques Appliqu\'ees (UMR 8050)
UPEM, F-77454, Marne-la-Vall\'ee, France}

\email{alain.pajor@@u-pem.fr}

\begin{abstract} We prove that there exists an absolute constant $C$ so that 
$$\mu(K)\ \le\ C\sqrt{p} \max_{\xi\in S^{n-1}} \mu(K\cap \xi^\bot)\ |K|^{1/n}$$
for any $p>2,$ any $n\in \N,$ any convex body $K$ that is the unit ball of an $n$-dimensional
subspace of $L_p,$ and any measure $\mu$ with non-negative even continuous density in $\R^n.$
Here $\xi^\bot$ is the central hyperplane perpendicular to a unit vector $\xi\in S^{n-1},$
and $|K|$ stands for volume. 
\end{abstract} 

\maketitle

\section{Introduction}  The slicing problem \cite{Bo1, Bo2, Ba1, MP}, a major open question  in convex geometry,
asks whether there exists a constant $C$ so that for any $n\in \N$ and any origin-symmetric convex body
$K$ in $\R^n,$
$$|K|^{\frac {n-1}n} \le C \max_{\xi\in S^{n-1}} |K\cap \xi^\bot|,$$
where $|K|$ stands for volume of proper dimension, and $\xi^\bot$ is the central hyperplane
in $\R^n$ perpendicular to a unit vector $\xi.$ The best-to-date
result $C\le O(n^{1/4})$ is due to Klartag \cite{Kl}, who
improved an earlier estimate of Bourgain \cite{Bo3}. The answer is affirmative
for unconditional convex bodies (as initially observed by Bourgain; see also \cite{MP, J2,
BN}), intersection bodies
\cite[Theorem 9.4.11]{G}, zonoids, duals of bodies with bounded volume ratio
\cite{MP}, the Schatten classes \cite{KMP}, $k$-intersection bodies \cite{KPY, K6};
see [BGVV] for more details.

The case of unit balls of finite dimensional subspaces of $L_p$ is of particular interest in this note.
It was shown by Ball \cite{Ba2} that the slicing problem has an affirmative answer for the unit
balls of finite dimensional subspaces 
of $L_p,\ 1\le p \le 2.$ Junge \cite{J1} extended this result to every $p\in (1,\infty),$ with the constant $C$
depending on $p$ and going to infinity when $p\to \infty.$ Milman \cite{M1} gave a different proof for
subspaces of $L_p,\ 2 < p <\infty,$ with the constant $C\le O(\sqrt{p}).$ Another proof of this estimate
can be found in \cite{KPY}. 

A generalization of the slicing problem to arbitrary measures was considered in \cite{K3,K4,K5,K6}. 
Does there exist a constant $C$ so that for every $n\in N,$
every origin-symmetric convex body $K$ in $\R^n,$ and every measure $\mu$ 
with non-negative even continuous density $f$ in $\R^n,$
\begin{equation}\label{mainproblem}
\mu(K)\ \le\ C \max_{\xi\in S^{n-1}} \mu(K\cap \xi^\bot)\ |K|^{1/n},
\end{equation}
where  $\mu(K)=\int_K f,$ and $\mu(K\cap \xi^\bot)=\int_{K\cap \xi^\bot} f ?$

Inequality (\ref{mainproblem}) was proved with an absolute constant $C$ for intersection bodies \cite{K3}
(by \cite{K2}, this includes the unit balls of subspaces of $L_p$ with $0<p\le 2)$,
unconditional bodies and duals of bodies with bounded volume ratio in \cite{K6},
for $k$-intersection bodies in \cite{K5}. For arbitrary origin-symmetric convex bodies, (\ref{mainproblem})
was proved in \cite{K4} with $C\le O(\sqrt{n}).$ A different proof of the latter estimate was recently given
in \cite{CGL}, where the symmetry condition was removed.

For the unit balls of subspaces of $L_p,\ p>2,$ (\ref{mainproblem}) was proved in \cite{K5}
with $C\le O(n^{1/2-1/p}).$ In this note we improve the estimate to $C\le O(\sqrt{p}),$ extending
Milman's result \cite{M1} to arbitrary measures in place of volume. In fact, we prove a more general
inequality
\begin{equation}\label{k-mainproblem}
\mu(K)\ \le\ (C\sqrt{p})^k \max_{H\in Gr_{n-k}} \mu(K\cap H)\ |K|^{k/n},
\end{equation}
where $1\le k < n,$ $Gr_{n-k}$ is the Grassmanian of $(n-k)$-dimensional subspaces
of $\R^n,$ $K$ is the unit ball of any $n$-dimensional subspace of $L_p,\ p>2,$
$\mu$ is a measure on $\R^n$ with even continuous density, and $C$ is a constant
independent of $p,n,k,K,\mu.$

The proof is a combination of two known results. Firstly, we use the reduction of the slicing
problem for measures to computing the outer volume ratio distance from a body to
the class of intersection bodies established in \cite{K6}; see Proposition \ref{lowdim}.
Note that outer volume ratio estimates have been applied to different cases of the original slicing problem
by Ball \cite{Ba2}, Junge \cite{J1}, and E.Milman \cite{M1}.
Secondly, we use an estimate for the outer volume ratio distance from the unit ball of a subspace of $L_p,\ p>2,$
to the class of origin-symmetric ellipsoids proved by E.Milman in \cite{M1}. This estimate also follows from results
of Davis, V.Milman and Tomczak-Jaegermann \cite{DMT}. We include a concentrated version
of the proof in Proposition \ref{ovr}.

\section{Slicing inequalities} \label{slicing}
We need several definitions and facts.
A closed bounded set $K$ in $\R^n$ is called a {\it star body}  if 
every straight line passing through the origin crosses the boundary of $K$ 
at exactly two points different from the origin, the origin is an interior point of $K,$
and the {\it Minkowski functional} 
of $K$ defined by 
$$\|x\|_K = \min\{a\ge 0:\ x\in aK\}$$
is a continuous function on $\R^n.$ 

The {\it radial function} of a star body $K$ is defined by
$$\rho_K(x) = \|x\|_K^{-1}, \qquad x\in \R^n,\ x\neq 0.$$
If $x\in S^{n-1}$ then $\rho_K(x)$ is the radius of $K$ in the
direction of $x.$

We use the polar formula for volume of a star body
\begin{equation}\label{polar}
|K|=\frac 1n \int_{S^{n-1}} \|\theta\|_K^{-n} d\theta.
\end{equation}

The class of intersection bodies was introduced by Lutwak \cite{L}.
Let $K, L$ be origin-symmetric star bodies in $\R^n.$ We say that $K$ is the 
intersection body of $L$ if the radius of $K$ in every direction is 
equal to the $(n-1)$-dimensional volume of the section of $L$ by the central
hyperplane orthogonal to this direction, i.e. for every $\xi\in S^{n-1},$
$$
\rho_K(\xi)= \|\xi\|_K^{-1} = |L\cap \xi^\bot|
$$
$$= \frac 1{n-1} \int_{S^{n-1}\cap \xi^\bot} \|\theta\|_L^{-n+1}d\theta=
\frac 1{n-1} R\left(\|\cdot\|_L^{-n+1}\right)(\xi),$$
where $R:C(S^{n-1})\to C(S^{n-1})$ is the {\it spherical Radon transform}
$$Rf(\xi)=\int_{S^{n-1}\cap \xi^\bot} f(x) dx,\qquad \forall f\in C(S^{n-1}).$$
All bodies $K$ that appear as intersection bodies of different star bodies
form {\it the class of intersection bodies of star bodies}. A more general class of {\it intersection bodies} 
is defined as follows. If $\mu$ is a finite Borel measure on $S^{n-1},$ then the spherical Radon transform
$R\mu$ of $\mu$ is defined as a functional on $C(S^{n-1})$ acting by
$$(R\mu, f)=(\mu, Rf)=\int_{S^{n-1}} Rf(x) d\mu(x),\qquad \forall f\in C(S^{n-1}).$$
A star body $K$ in $\R^n$ is called an {\it intersection body} if $\|\cdot\|_K^{-1}=R\mu$
for some measure $\mu,$ as functionals on $C(S^{n-1}),$  i.e.
$$\int_{S^{n-1}} \|x\|_K^{-1} f(x) dx = \int_{S^{n-1}} Rf(x)d\mu(x),\qquad \forall f\in C(S^{n-1}).$$
Intersection bodies played a crucial role in the solution of the
Busemann-Petty problem and its generalizations; see \cite[Chapter 5]{K1}.

A generalization of the concept of an intersection body
was introduced by Zhang \cite{Z} 
in connection with the lower dimensional Busemann-Petty problem.
For $1\le k \le n-1,$  the {\it $(n-k)$-dimensional spherical Radon transform} 
$R_{n-k}:C(S^{n-1})\to C(Gr_{n-k})$  
is a linear operator defined by
$$R_{n-k}g (H)=\int_{S^{n-1}\cap H} g(x)\ dx,\quad \forall  H\in Gr_{n-k}$$
for every function $g\in C(S^{n-1}).$
\smallbreak
We say that
an origin symmetric star body $K$ in $\R^n$ is a {\it generalized $k$-intersection body}, 
and write $K\in {\cal{BP}}_k^n,$  if there exists a finite Borel non-negative measure $\mu$
on $Gr_{n-k}$ so that for every $g\in C(S^{n-1})$
\begin{equation}\label{genint}
\int_{S^{n-1}} \|x\|_K^{-k} g(x)\ dx=\int_{Gr_{n-k}} R_{n-k}g(H)\ d\mu(H).
\end{equation}
When $k=1$ we get the class of intersection bodies.
It was proved by Goodey and Weil \cite{GW} for $k=1$ and by Grinberg and Zhang \cite[Lemma 6.1]{GZ}
for arbitrary $k$ (see also \cite{M2} for a different proof) that the class ${\cal{BP}}_k^n$ is the closure in the radial metric
of $k$-radial sums of origin-symmetric ellipsoids. In particular, the classes ${\cal{BP}}_k^n$ contain all origin-symmetric
ellipsoids in $\R^n$ and are invariant with respect to linear transformations. 
Recall that the $k$-radial sum $K+_kL$ of star bodies $K$ and $L$ is defined by
$$\rho_{K+_kL}^{k}= \rho_K^{k} + \rho_L^k.$$

For a convex body $K$ in $\R^n$ and $1\le k <n,$ denote by 
$${\rm {o.v.r.}}(K,{\cal{BP}}_k^n) = \inf \left\{ \left( \frac {|C|}{|K|}\right)^{1/n}:\ K\subset C,\ C\in {\cal{BP}}_k^n \right\}$$
the outer volume ratio distance from a body $K$ to the class ${\cal{BP}}_k^n.$

Let $B_2^n$ be the unit Euclidean ball in $\R^n,$ let $|\cdot|_2$ be the Euclidean norm in $\R^n,$
and let $\sigma$ be the uniform probability measure on the sphere $S^{n-1}$ in $\R^n.$
For every $x\in\R^n$, let $x_1$ be the first coordinate of $x$. We use the fact that for every $p>-1$
\begin{equation} \label{moment}
\int_{S^{n-1}} |x_1|^p d\sigma(x) = \frac{\Gamma(\frac{p+1}2) \Gamma(\frac n2)}
{\sqrt{\pi} \Gamma(\frac{n+p}2)};
\end{equation}
see for example \cite[Lemma 3.12]{K1}, where one has to divide by  
$|S^{n-1}|=  2\pi^{(n-1)/2}/\Gamma(\frac n2),$ because the measure $\sigma$ on the sphere
is normalized.
 
In \cite{K6}, the slicing problem for arbitrary measures was reduced to estimating the outer volume
ratio distance from a convex body to the classes   ${\cal{BP}}_k^n$, as follows.
\begin{pr} \label{lowdim} For any $n\in \N,\ 1\le k <n,$ any origin-symmetric star body $K$ in $\R^n,$ 
and any measure $\mu$ with even continuous density on $K,$
$$\mu(K)\le \left({\rm{ o.v.r.}}(K,{{\cal{BP}}_k^n})\right)^k\ \frac n{n-k} 
c_{n,k} \max_{H\in Gr_{n-k}} \mu(K\cap H)\ |K|^{k/n},$$
where $c_{n,k}=|B_2^n|^{(n-k)/n}/|B_2^{n-k}| \in (e^{-k/2},1).$
\end{pr}

It appears that for the unit balls of subspaces of $L_p,\ p>2$ the outer volume ration distance
to the classes of intersection bodies does not depend on the dimension. As mentioned in the
introduction, the following estimate was proved in \cite{M1} and also follows from results of \cite{DMT}.
We present a short version of the proof.

\begin{pr} \label{ovr} Let $p>2,\ n\in \N,\ 1\le k <n,$ and let $K$ be the unit ball of 
an $n$-dimensional subspace of $L_p.$ Then
$${\rm{ o.v.r.}}(K,{{\cal{BP}}_k^n})\le C\sqrt{p},$$
where $C$ is an absolute constant.
\end{pr}

\pf Since the classes ${{\cal{BP}}_k^n}$ are invariant under linear transformations, we can assume that 
$K$ is in the Lewis position. By a result of Lewis in the form of \cite[Theorem 8.2]{LYZ}, 
this means that there exists
a measure $\nu$ on the sphere so that for every $x\in \R^n$
$$\|x\|_K^p= \int_{S^{n-1}} |(x,u)|^p d\nu(u),$$
and
$$|x|_2^2=\int_{S^{n-1}} |(x,u)|^2 d\nu(u).$$
Also, by the same result of Lewis \cite{Le}, $K\subset n^{1/2-1/p}B_2^n.$ 

Let us estimate volume of $K$ from below. By the Fubini theorem, formula (\ref{moment}) and Stirling's
formula, we get
$$\int_{S^{n-1}} \|x\|_K^p d\sigma(x) =  \int_{S^{n-1}} \int_{S^{n-1}} |(x,u)|^p d\sigma(x) d\nu(u)$$
$$=\int_{S^{n-1}} |x_1|^p d\sigma(x) \int_{S^{n-1}} d\nu(u) \le \left(\frac{Cp}{n+p}\right)^{p/2} \int_{S^{n-1}} d\nu(u).$$
Now
$$\frac{Cp}{n+p}\left(\int_{S^{n-1}} d\nu(u)\right)^{2/p}\ge \left(\int_{S^{n-1}} \|x\|_K^p d\sigma(x) \right)^{2/p}$$
$$\ge \left(\int_{S^{n-1}} \|x\|_K^{-n} d\sigma(x) \right)^{-2/n} = \left(\frac{|K|}{|B_2^n|}\right)^{-2/n}\sim \frac1n |K|^{-2/n},$$
because $|B_2^n|\sim n^{-1/2}.$
On the other hand,
$$1= \int_{S^{n-1}} |x|_2^2 d\sigma(x) = \int_{S^{n-1}} \int_{S^{n-1}} (x,u)^2 d\nu(u) d\sigma(x)$$
$$= \int_{S^{n-1}} \int_{S^{n-1}} |x_1|^2 d\sigma(x) d\nu(u)= \frac 1n \int_{S^{n-1}} d\nu(u),$$
so
$$\frac{Cp}{n+p} n^{2/p}\ge \frac 1n |K|^{-2/n},$$
and
$$|K|^{1/n}\ge c n^{-1/p} \sqrt{\frac{n+p}{np}} \ge \frac {cn^{1/2-1/p}}{\sqrt{p}} |B_2^n|^{1/n}.$$

Finally, since $K\subset n^{1/2-1/p}B_2^n,$ and $B_2^n\in {\cal{BP}}_k^n$ for every $k,$ we have
$${\rm{ o.v.r.}}(K,{{\cal{BP}}_k^n}) \le \left(\frac{|n^{1/2-1/p}B_2^n|}{|K|}\right)^{1/n}\le C\sqrt{p},$$
where $C$ is an absolute constant.
\bigbreak

We now formulate the main result of this note. 
\begin{co} There exists a constant $C$ so that for any $p>2,\ n\in \N,\ 1\le k<n,$ any convex body $K$
that is the unit ball of an $n$-dimensional subspace of $L_p,$ and any measure $\mu$ with non-negative even
continuous density in $\R^n,$
$$\mu(K)\ \le\ (C\sqrt{p})^k \max_{H\in Gr_{n-k}} \mu(K\cap H)\ |K|^{k/n}.$$
\end{co}

\pf Combine Proposition \ref{lowdim} with Proposition \ref{ovr}. Note that
$\frac n{n-k}\in (1, e^k),$ and $c_{n,k}\in (e^{-k/2},1),$ so these constants
can be incorporated in the constant $C.$
\endpf

\end{document}